\documentclass{elsart}
\usepackage{amsmath,amsfonts,amssymb}

\newcommand\Rn{{\mathbb R}^n}
\newcommand{\irm}{{\rm i}}

\renewcommand{\S}{\mathcal{S}}

\newcommand{\mb}[1]{\ensuremath{\mathbb{#1}}}
\newcommand{\N}{\mb{N}}
\newcommand{\R}{\mb{R}}

\newcommand{\lara}[1]{\langle #1 \rangle}
\newcommand{\esp}{\mathrm{e}}
\newcommand{\eps}{\varepsilon}
\newcommand{\beq}{\begin{equation}}
\newcommand{\eeq}{\end{equation}}

\begin{document}

\begin{frontmatter}
\title{On the well-posedness of weakly hyperbolic equations with time dependent coefficients}
\author[IC]{Claudia Garetto\thanksref{th:grant1}}
\author[IC]{Michael Ruzhansky\thanksref{th:grant2}}
\thanks[th:grant1]{The first
 author was supported by
the Imperial College Junior Research Fellowship.}
 \thanks[th:grant2]{The second author was supported by the
EPSRC Leadership Fellowship EP/G007233/1.}
\address[IC]{Department of Mathematics,\\
Imperial College London,\\
180 Queen's Gate, London SW7 2AZ, UK}

\begin{abstract}
In this paper we analyse the Gevrey well-posedness of the Cauchy problem
for weakly hyperbolic equations of general form with time-dependent coefficients.
The results involve the order of lower order terms and the number of multiple
roots. We also derive the corresponding well-posedness results in the space
of Gevrey Beurling ultradistributions.
\end{abstract}
\end{frontmatter}

\section{Introduction}

In this paper we study the well-posedness for weakly hyperbolic equations of higher
orders of general form with time-dependent coefficients. Namely, we consider the Cauchy problem
\beq
\label{CP}
\left\{
\begin{array}{cc}
D^m_t u=\sum_{j=0}^{m-1} A_{m-j}(t,D_x)D_t^j u+f(t,x),&\quad (t,x)\in[0,T]\times\R^n,\\
D^{k-1}_t u(0,x)=g_{k}(x),&\quad k=1,...,m,
\end{array}
\right.
\eeq
where each $A_{m-j}(t,D_x)$ is a differential operator of order $m-j$ with continuous coefficients only depending on $t$.
As usual, $D_t=\frac{1}{\irm}\partial_t$ and $D_x=\frac{1}{\irm}\partial_x$.
More precisely, we can write equation \eqref{CP} as
\beq
\label{CP1}
D^m_t u=\sum_{j=0}^{m-1}\sum_{|\gamma|=m-j} a_{m-j,\gamma}(t)D_x^\gamma
D_t^j u+
\sum_{|\gamma|+j\leq l} a_{m-j,\gamma}(t)D_x^\gamma
D_t^j u+f(t,x),
\eeq
where $l$ is the order of lower order terms, $0\leq l\leq m-1$.
Concerning the lower order terms, throughout the paper we will only assume that
$a_{m-j,\gamma}(t)\in C[0,T]$ for $|\gamma|+j\leq l$, and that
$f\in C([0,T];G^s(\Rn))$ is continuous in $t$ and Gevrey in $x$ of order
$s$ appearing in the formulation of the theorems below.

Weakly hyperbolic equations \eqref{CP}, \eqref{CP1} and their special cases have been extensively considered
in the literature, see e.g. \cite{B,  ColKi:02, ColKi:02-2, DK, I, KS}, to mention only very few, and references therein.

Let $A_{(m-j)}$ denote the principal part of the operator $A_{m-j}$ and let $\tau_k(t,\xi)$, $k=1,...,m$, be the roots of the characteristic equation
\[
\tau^m=\sum_{j=0}^{m-1}A_{(m-j)}(t,\xi)\tau^j
\equiv\sum_{j=0}^{m-1}\sum_{|\gamma|=m-j}a_{m-j,\gamma}(t)\xi^\gamma\tau^j.
\]
We will analyse the following two cases:

\medskip
%\begin{itemize}
%\item[
\boxed{\bf Case 1} we assume that
%\begin{itemize}
%\item[{\,}]
\\[0.3cm] the roots $\tau_k(t,\xi)$, $k=1,...,m$, are real-valued and of H\"older class $C^\alpha$, $0<\alpha\le 1$
with respect to $t$; for any $t\in[0,T]$ they either coincide or are all distinct.
%\end{itemize}
%\item[

\medskip
\boxed{\bf Case 2} there exists $r=2,...,m-1$ such that
\begin{itemize}
\item[(i)] the roots $\tau_k(t,\xi)$, $k=1,...,r$, are real-valued, of class $C^\alpha$, $0<\alpha\le 1$ with respect to $t$ and either coincide or are all distinct;
\item[(ii)] the roots $\tau_k(t,\xi)$, $k=r+1,...,m$, are real-valued, of class $C^\beta$, $0<\beta\le 1$ with respect to $t$ and are all distinct.
\end{itemize}
%\end{itemize}
Before we proceed we note that in the case $\alpha=1$ or $\beta=1$, it is enough
to assume Lipschitz regularity for the corresponding roots. This includes
the case of weakly hyperbolic equations with smooth coefficients in which
case the roots are Lipschitz by Bronshtein's theorem.

In the next section we give the Gevrey well-posedness results for
the Cauchy problem \eqref{CP} in Case 1 and Case 2, as well as in the
strictly hyperbolic case formulated below in Case 3. In summary, our
results will apply to all dimensions and will improve the known Gevrey
indices in different settings. First we describe what is known for this
problem.

Cauchy problems of such type
have been studied in the Gevrey framework by Colombini and Kinoshita in \cite{ColKi:02} but only in the one dimensional case, i.e., $x\in\R$, and with $f\equiv 0$. In the present paper, we extend the result of \cite{ColKi:02} to any
dimension $n\geq 1$, as well as improve the indices for the Gevrey well-posedness
(see Remarks \ref{REM:case1} and \ref{REM:case2}).
The idea of the proof in \cite{ColKi:02} is to reduce the Cauchy problem
 \eqref{CP} to a differential system keeping track of all the derivatives of the
solution $u$. The new unknown function contains also the lower order derivatives of $u$ and
thus the size of the resulting system is much higher than $m$. Technically, it makes it
hard to extend this method to higher dimensions. In this paper we use the
pseudo-differential techniques of the reduction of \eqref{CP} to the system. This allows
us to keep the size of the system to be equal to $m$ and works equally well in all dimensions.
The subsequent estimates can be then improved for several terms in the proof of the
energy inequality.
Here,  we also give results for inhomogeneous equations as well
as discuss the well-posedness of the problem \eqref{CP}
in the spaces of ultradistributions.

More generally, in dimensions $n\geq 1$, there are a number of results available
concerning the problem \eqref{CP}. It was known since a long time
(see Ivrii \cite{I} and references therein) that the Cauchy problem for
any hyperbolic equation with sufficiently smooth coefficients is well-posed
in Gevrey classes $G^s$ with $1\leq s< s_m$ for some $s_m>1$.
Subsequently, it was shown by Bronshtein \cite{B}
that the equation \eqref{CP} with characteristics of multiplicity $r\leq m$,
with coefficients $C^\infty$ in $t$ (and also allowing $G^s$ in $x$), is
well-posed in $G^s$ for $1\leq s<1+\frac{1}{r-1}$.
This bound is in general sharp but can be improved in particular cases,
such as, for example, Case 1 in Theorem \ref{THM:case1}, allowing
lower regularity on the coefficients and taking into account the degree of
lower order terms.
Under the smoothness assumptions on the coefficients, we have
$\alpha=\beta=1$ in our assumptions, so that
the index $\frac{1}{r-1}$ corresponds to $\frac{\beta}{r-\beta}$ in
Theorem \ref{THM:case2}.

When $m=2$, $l=1$ and $r=2$, Colombini, De Giorgi, Jannelli and Spagnolo
(see \cite{CDS, CJS}) considered equations \eqref{CP}, \eqref{CP1} with
$a_{1,1}(t)\equiv 0$ and $a_{0,2}(t)\in C^\delta[0,T]$, $\delta> 0$.
They showed that the Cauchy problem
\eqref{CP} is well-posed in $G^s$ provided that $1\leq s<1+\frac{\delta}{2}$.
In our setting this is covered by the conditions of Case 1 with
$\alpha=\frac{\delta}{2}$, so that the result above is included in
Theorem \ref{THM:case1} giving the range
$1\leq s<1+\alpha$. They also considered the case of $r=1$ when
they proved the well-posedness in $G^s$ for
$1\leq s<1+\frac{\delta}{1-\delta}$.
In our setting this falls under the assumptions of Case 2
with $\alpha=\frac{\delta}{2}$ and $\beta=\delta$,
so that their result is included in Theorem \ref{THM:case2}
with the same range for $s$.

In \cite{KS}, Kinoshita and Spagnolo considered the Cauchy problem \eqref{CP}
for operators with homogeneous symbols in one-dimension, i.e.
assuming that $n=1$ and $a_{m-j,\gamma}(t)\equiv 0$ for
$\gamma+j<m$ in \eqref{CP1}. Among other results for such equations, they
showed that if $a_{m-j,\gamma}(t)\in C^2[0,T]$, $\gamma+j=m$,
and the characteristic roots satisfy
\beq\label{KS}
\tau_i(t)^2+\tau_j(t)^2\leq M(\tau_i(t)-\tau_j(t))^2 \quad\textrm{ for } i\not=j,
\eeq
then
the Cauchy problem \eqref{CP} is well-posed in the Gevrey space $G^s$
provided that $1\leq s<1+\frac{1}{m-1}$. In our setting,
the condition $a_{m-j,\gamma}(t)\in C^2[0,T]$ corresponds to
$\alpha=\frac{2}{r}$ and $\beta=1$. Thus,
Theorem \ref{THM:case2} implies the well-posedness in the Gevrey space
$G^s$ for $1\leq s<1+\min\{\frac{2}{r},\frac{1}{r-1}\}=1+\frac{1}{r-1}$,
provided that the equation has multiplicities ($2\leq r\leq m$).
In this sense the result of Theorem \ref{THM:case2} improves the
$C^2$-coefficients result of \cite{KS}, also allowing any $n\geq 1$ and lower order terms,
as well as removing the assumption \eqref{KS} on the roots. {We note that
condition \eqref{KS} has been considered earlier in Colombini and Orr\`u in \cite{COr} to prove
$C^\infty$ well-posedness in case of analytic coefficients.} Certain improvements
have been also observed by Jannelli in \cite{Jan:09}.

We also present the corresponding results for the well-posedness in classes of Gevrey
ultradistributions. It is by now well-known that the Cauchy problems for weakly
hyperbolic equations even with smooth coefficients do not have to be in general well-posed in
the space ${\mathcal D}'(\Rn)$ of distributions, see e.g. Colombini, Jannelli and
Spagnolo \cite{CJS} and Colombini, Spagnolo \cite{CS}. In the subsequent paper
we will analyse the propagation of singularities for weakly hyperbolic equations,
and for such purpose it is necessary to have a framework in which the Cauchy
problem would be well-posed. In fact, such a well-posedness result follows
directly from the energy estimates that we will establish in the proofs of
Theorem \ref{THM:case1} and Theorem \ref{THM:case2}. However, there
is still one subtle matter of the definition of the corresponding space of Gevrey
ultradistributions. Namely, we will show that one has to take the Beurling Gevrey
ultradistributions rather than the Roumieu Gevrey ultradistributions to achieve
such results. In general, in the absence of energy inequalities certain conclusions
in spaces containing Schwartz distributions are also possible, but such questions
will be treated elsewhere.

Furthermore, we complement the weakly hyperbolic analysis by giving the
results for strictly hyperbolic equations with coefficients of low regularity.
This corresponds to Case 2 above when we take $r=1$. As the equation is
strictly hyperbolic, we do not have to distinguish between regularities
of simple and multiple roots, so that we can take $\alpha=\beta$ in
this case. To summarise, we consider
%\begin{itemize}
%\item[

\medskip
\boxed{\bf Case 3} we assume that
%\begin{itemize}
%\item[{\,}]
\\[0.3cm]
the roots $\tau_k(t,\xi)$, $k=1,...,m$, are real-valued and of H\"older class $C^\beta$, $0<\beta\le 1$
with respect to $t$; for any $t\in[0,T]$ they are all distinct.
%\end{itemize}
%\end{itemize}
%

\medskip
The proof of the corresponding statements will follow by taking the
proof of Case 2 and putting $r=1$ and $\alpha=\beta$ at the end.

Finally we note that if the operator in \eqref{CP} is strictly hyperbolic
and coefficients are more regular,
much more is known. For a detailed analysis of large-time asymptotics
of equations \eqref{CP}, \eqref{CP1}, for constant coefficients,
we refer to Ruzhansky and Smith \cite{RS}. Here we
note that althought \eqref{CP} may be strictly hyperbolic, multiplicities
of the full equation (together with lower order terms)
may still occur for small frequencies due to the presence of low
order terms. Equations with $C^1$-regularity of the coefficients
with respect to time
have been treated in Matsuyama and Ruzhansky \cite{MR}, while systems with oscillations and
more regularity have been analysed in Ruzhansky and Wirth \cite{RW}.

In the sequel, we denote $\langle\xi\rangle=(1+|\xi|^{2})^{1/2}.$

The authors thank Professor T. Kinoshita for useful discussions,
and the Daiwa
foundation for support.

\numberwithin{equation}{section}

\section{Main results}

From the fact that each $A_{(m-j)}(t,\xi)$ is a polynomial homogeneous of degree $m-j$ in $\xi$ it follows that the roots $\tau_k(t,\xi)$ are positively homogeneous of degree $1$ in $\xi$. Combining this fact with the H\"older regularity it
follows in Case 1 that there exists a constant $c>0$ such that
\beq
\label{hyp_roots}
|\tau_k(t,\xi)-\tau_k(s,\xi)|\le c|\xi||t-s|^\alpha
\eeq
for $k=1,\ldots,m$, for all $\xi\neq 0$ and $t,s\in[0,T]$.

Throughout the paper, without loss of generality, by relabeling the roots,
we can always arrange that they are ordered, so that we will assume that
\beq
\label{hyp_roots_2}
\tau_1(t,\xi)\le\tau_2(t,\xi)\le\cdots\le\tau_m(t,\xi),
\eeq
for all $t$ \and $\xi$.
The index for the H\"older regularity is preserved under such a relabelling.
More precisely, in Case 2 we have that \eqref{hyp_roots} is true with exponent $\beta$ when $r+1\le k\le m$ and
\begin{multline}
\label{hyp_roots_3}
\tau_1(t,\xi)\le\tau_2(t,\xi)\le\cdots\le\tau_r(t,\xi)<\tau_{r+1}(t,\xi)
\\
<\tau_{r+2}(t,\xi)<\cdots<\tau_m(t,\xi),
\end{multline}
for all $t$ and $\xi\neq 0$. From the homogeneity in $\xi$ we also have that
\beq
\label{strict_hyp}
\tau_{k}(t,\xi)-\tau_{k-1}(t,\xi)\ge c_0|\xi|
\eeq
for some constant $c_0>0$ and all $k=r+1,...,m$, uniformly in $t\in[0,T]$ and $\xi\neq 0$.
Throughout the paper we also
assume that the roots which coincide have the following uniform property: there exists $c>0$ such that \beq
\label{hyp_coincide}
|\tau_i(t,\xi)-\tau_j(t,\xi)|\le c|\tau_k(t,\xi)-\tau_{k-1}(t,\xi)|
\eeq
for all $1\le i,j,k\le r$, for all $t\in[0,T]$ and $\xi\in\R^n$.
We note that although
condition \eqref{hyp_coincide} was not explicitly stated in
\cite{ColKi:02}, it is required for their
proof also in the case $n=1$
(\cite{Kp}).

We first formulate the results in Gevrey spaces.
Throughout the formulations, in inequalities for indices, we will adopt the convention
that $\frac{1}{0}=+\infty$. We briefly recall the definition of the space
$G^s(\Rn)$, the space of (Roumieu) Gevrey functions.
We denote $\N_0=\N\cup\{0\}$.

\begin{defn}\label{def_gevrey}
Let $s\geq 1$.
We say that $f\in C^\infty(\R^n)$ belongs to the Gevrey class $G^s(\R^n)$ if for every compact set $K\subset\R^n$ there
exists a constant $C>0$ such that for all $\alpha\in\N_0^n$ we have the estimate
\[
\sup_{x\in K}|\partial^\alpha f(x)|\le C^{|\alpha|+1}(\alpha!)^s.
\]
\end{defn}
We recall that $G^1(\R^n)$ is the space of analytic functions and that
$G^s(\R^n)\subseteq G^\sigma(\R^n)$ if $s\le \sigma$.
For  $s>1$, let
$G^s_0(\R^n)$ be the space of compactly supported Gevrey functions of order $s$. In the paper we make use of the following Fourier characterisation (see \cite[Theorem 1.6.1]{Rodino:93}), where $\lara{\xi}=(1+|\xi|^2)^{\frac{1}{2}}$.
\begin{prop}
\label{prop_Fourier}
\leavevmode
\begin{itemize}
\item[(i)] Let $u\in G^s_0(\R^n)$. Then, there exist constants $c>0$ and $\delta>0$ such that
\beq
\label{fou_gey_1}
|\widehat{u}(\xi)|\le c\,\esp^{-\delta\lara{\xi}^{\frac{1}{s}}}
\eeq
for all $\xi\in\R^n$.
\item[(ii)] Let $u\in\S'(\R^n)$. If there exist constants $c>0$ and $\delta>0$ such that \eqref{fou_gey_1} holds then $u\in G^s(\R^n)$.
    \end{itemize}
\end{prop}

We now formulate the result for Case 1.

\begin{thm}\label{THM:case1}
Let $T>0$ and $0\leq l\leq m-1$. Assume the conditions of Case 1.
Then for any $g_k(x)\in G^s(\Rn)$ ($k=1,\ldots,m$),
the Cauchy problem \eqref{CP} has a unique global solution
$u\in C^m([0,T]; G^s(\Rn))$, provided that
\begin{equation}\label{EQ:case1-s}
1\leq s<1+\min\left\{\alpha,\frac{m-l}{l}\right\}.
\end{equation}
\end{thm}

\begin{rem}\label{REM:case1}
In \cite{ColKi:02}, the authors proved that in one dimension $n=1$ and
for $f\equiv 0$, the well-posedness in
Theorem \ref{THM:case1} holds provided that
$1\leq s < 1+\min\left\{\alpha,\frac{m-l}{m-1}\right\}$.
Theorem  \ref{THM:case1} improves this result by increasing the second factor
under the minimum, as well as gives the result for any dimensions
and non-zero $f$.
\end{rem}

If we observe that $\alpha\leq \frac{m-l}{l}$ is equivalent to
$l\leq \frac{m}{\alpha+1}$, we get
\begin{cor}\label{lot}
Under conditions of Case 1, if the order of lower order terms satisfies
$l\leq \frac{m}{2}$, then the well-posedness in  Theorem \ref{THM:case1}
holds for $1\leq s<1+\alpha$. More precisely, if we assume that $l\leq \frac{m}{\alpha+1}$, then the
well-posedness in Theorem \ref{THM:case1} holds provided that
$1\leq s<1+\alpha$.
\end{cor}

Under assumptions of Case 2, if there are simple roots, sometimes the index in \eqref{EQ:case1-s}
can be improved. However, this should not be generally expected as a multiplication of
a weakly hyperbolic polynomial by a strictly hyperbolic one should not, in general,
improve the well-posedness of the Cauchy problem.

\begin{thm}\label{THM:case2}
Let $T>0$, $2\leq r\leq m-1$ and $0\leq l\leq m-1$. Assume the conditions of Case 2.
Then for any $g_k(x)\in G^s(\Rn)$ ($k=1,\ldots,m$),
the Cauchy problem \eqref{CP} has a unique global solution
$u\in C^m([0,T]; G^s(\Rn))$, provided that
\begin{equation}\label{EQ:case2-s}
1\leq s<1+\min\left\{\alpha,\frac{\beta}{r-\beta}\right\}.
\end{equation}
\end{thm}

\begin{rem}\label{REM:case2}
In \cite{ColKi:02}, the authors proved that in one dimension $n=1$ and
for $f\equiv 0$, the well-posedness in
Theorem \ref{THM:case2} holds provided that
$1\leq s < 1+\min\left\{\alpha,\frac{\beta}{r-\beta},\frac{m-l}{r-1}\right\}$.
Theorem  \ref{THM:case2} improves this result by removing the last term
under the minimum, as well as applies to all dimensions and non-zero $f$.
\end{rem}

We now give a remark about the strictly hyperbolic equations covered by
Case 3. We recall that in this case we take $\alpha=\beta$.

\begin{rem}\label{REM:r1statemement}
Under the conditions of Case 3, the conclusion of
Theorem  \ref{THM:case2} holds provided that $$1\le s<1+\frac{\beta}{1-\beta}.$$
See Remark \ref{REM:r1} for the proof.
\end{rem}

The result of Theorem \ref{THM:case2} is better than that in Theorem \ref{THM:case1}
if there are few multiple roots, or if the order of lower order terms is
sufficiently high. In particular and more precisely,
it can be easily checked  that $\frac{\beta}{r-\beta}\geq \frac{m-l}{l}$ if
$r\leq \frac{\beta m}{m-l}$ (where $r$ is the number of multiple roots),
or if $l\geq \frac{(r-\beta)m}{r}$ (where $l$ is the order of lower order terms).

It is interesting to observe the implications of Theorem \ref{THM:case2} for equations with at most double
roots ($r=2$) in the Cauchy problem \eqref{CP}, \eqref{CP1}, where the coefficients
$a_{m-j,\gamma}$ belong to $C^\delta$ with $0<\delta\le 1$ and $|\gamma|+j=m$. In this case we have
$\alpha=\frac{\delta}{2}$ and $\beta=\delta$, and
since $\frac{\delta}{2}<\frac{\delta}{2-\delta}$,
we obtain

\begin{cor}
Assume that in Case 2, we have $r=2$ (i.e. double roots) and that
for $|\gamma|+j=m$ the coefficients satisfy $a_{m-j,\gamma}\in C^\delta[0,T]$, $0< \delta\leq 1$.
Then the Cauchy problem
\eqref{CP} is well-posed in $C^m([0,T]; G^s(\Rn))$ for $1\leq s<1+\frac{\delta}{2}$.
\end{cor}

Finally we observe that all the arguments in the proofs remain valid if the
equation \eqref{CP} is pseudo-differential in the $x$-variable:

\begin{rem}\label{REM:local}
The results of Theorems \ref{THM:case1} and \ref{THM:case2} apply in the same
way for operators in \eqref{CP} that are classical pseudo-differential in $D_x$, if we take
$g_k\in G^s_0(\Rn)$, $k=1,\ldots,m$, to be compactly supported.

Also, if the equality in \eqref{EQ:case1-s} or \eqref{EQ:case2-s} is attained,
the local well-posedness in Theorems \ref{THM:case1} and \ref{THM:case2}
and in the first part of this remark still holds.
\end{rem}

Before proceeding with the ultradistributions and with the
proof of the theorems above we give some examples. For more examples in
one dimension we refer to \cite{ColKi:02}, with the corresponding
improvement in indices given by Remarks \ref{REM:case1} and
\ref{REM:case2}.

\medskip
{\bf Example 1.} Let us consider the equation
\[
D_t^3u=-a(t)D_t \Delta_x  u+L(t,D_x,D_t)u,
\]
where $a(t)\ge 0$ belongs to $C^{2\alpha}([0,T])$, $\Delta_x=\partial_{x_1}^2+...+\partial_{x_n}^2$ and $L$ is a differential operator of order $l\le 2$. The corresponding principal symbol is $\tau^3-a(t)|\xi|^2\tau$ with roots
$\tau_1=-\sqrt{a(t)}|\xi|$, $\tau_2=0$ and $\tau_3=\sqrt{a(t)}|\xi|$. According to Theorem \ref{THM:case1} given initial data in $G^s(\R^n)$ the corresponding Cauchy problem has a unique solution $u\in C^3([0,T]; G^s(\Rn))$ with
\[
1\le s <1+\min\left\{\alpha,\frac{3-l}{l}\right\}.
\]
Note that the same well-posedness result holds for
\[
D_t^3u=\sum_{i=1}^n b_i(t)D_{x_i}D_t^2+L(t,D_x,D_t)u,
\]
when we assume that the coefficients $b_i$ are real-valued of class $C^{\alpha}$ and the multiplicity
is at a point $t_0\in[0,T]$ such that $b_i(t_0)=0$ for all $i=1,...,n$. We can apply Theorem \ref{THM:case1} and as
an example of the reordering of  the roots in the proof, we relabel the roots of the characteristic polynomial $\tau^3-\sum_i b_i(t)\xi_i\tau^2$ as
\[
\tau_1(t,\xi)=\min\left\{\sum_i b_i(t)\xi_i,0\right\},\, \tau_2=0,\, \tau_3(t,\xi)=\max\left\{\sum_i b_i(t)\xi_i,0\right\}.
\]

{\bf Example 2.} We study the Cauchy problem
\[
D_t^4 u=-(a(t)+b(t))D_t^2\Delta u-a(t)b(t)\Delta^2u,\qquad D_t^j u(0,x)=g_j(x),\, j=0,1,2,3,
\]
where we take $a\in C^{2\alpha}[0,T]$, $b\in C^{\beta}[0,T]$ with $a(t)\ge 0$ and $b(t)-a(t)\ge\delta>0$. The roots of the characteristic polynomial are $\tau_1(t,\xi)=-\sqrt{a(t)}|\xi|$, $\tau_2(t,\xi)=+\sqrt{a(t)}|\xi|$, $\tau_3(t,\xi)=-\sqrt{b(t)}|\xi|$ and $\tau_4(t,\xi)=+\sqrt{b(t)}|\xi|$. Hence, $r=2$ and from Theorem \ref{THM:case2} we have well-posedness in $C^4([0,T]; G^s(\Rn))$ with
\[
1\leq s<1+\min\left\{\alpha,\frac{\beta}{2-\beta}\right\}.
\]
Equations of this type were considered by Colombini and Kinoshita
in \cite{ColKi:02-2}, where the well-posedness was proved for
$1\leq s<1+\min\{\alpha,\frac{\beta}{2}\}.$ Thus, Theorem \ref{THM:case2} gives an improvement
of this result since $\frac{\beta}{2-\beta}\geq \frac{\beta}{2}$.
It also extends the one-dimensional version of this equation
considered in \cite[Example 3]{ColKi:02}.

Before stating the ultradistributional versions of Theorems \ref{THM:case1} and \ref{THM:case2} we recall a few more facts concerning Gevrey classes and ultradistributions. For more details see Komatsu \cite{K}, or
Rodino \cite[Section 1.5]{Rodino:93} for a partial treatment.
We first recall the Beurling Gevrey functions.

\begin{defn}\label{def_gevrey2}
Let $s\geq 1$.
We say that $f\in C^\infty(\R^n)$ belongs to the Beurling Gevrey class $G^{(s)}(\R^n)$ if for every compact set $K\subset\R^n$
and for every constant $A>0$ there exists a constant $C_{A,K}>0$ such that for all $\alpha\in\N_0^n$ we have the estimate
\[
\sup_{x\in K}|\partial^\alpha f(x)|\le C_{A,K} A^{|\alpha|}(\alpha!)^s.
\]
\end{defn}
Analogously to Proposition \ref{prop_Fourier}, we have the following Fourier characterisation, where $G^{(s)}_0(\R^n)$ denotes the space of compactly supported Beurling Gevrey functions.
\begin{prop}
\label{prop_Fourie2r}
\leavevmode
\begin{itemize}
\item[(i)] Let $u\in G^{(s)}_0(\R^n)$. Then, for any $\delta>0$ there
exists $C_\delta>0$ such that
\beq
\label{fou_gey_2}
|\widehat{u}(\xi)|\le C_\delta\,\esp^{-\delta\lara{\xi}^{\frac{1}{s}}}
\eeq
for all $\xi\in\Rn$.
\item[(ii)] Let $u\in\S'(\R^n)$. If for any $\delta>0$ there
exists $C_\delta>0$ such that \eqref{fou_gey_2} holds then $u\in G^{(s)}(\R^n)$.
\end{itemize}
\end{prop}
For $s>1$, the spaces $G^s_0(\Rn)$ and $G^{(s)}_0(\Rn)$ of compactly supported functions
can be equipped with natural seminormed topologies, and by $\mathcal{D}'_s(\R^n)$ and
$\mathcal{D}'_{(s)}(\R^n)$ we denote the spaces of linear continuous
functionals on them, respectively. We use the expressions Gevrey Roumieu ultradistributions
and Gevrey Beurling ultradistributions for the elements of $\mathcal{D}'_s(\R^n)$ and
$\mathcal{D}'_{(s)}(\R^n)$, respectively.
Let $\mathcal{E}'_s(\R^n)$ and $\mathcal{E}'_{(s)}(\R^n)$ be the topological duals of $G^s(\Rn)$
and $G^{(s)}(\Rn)$, respectively.
By duality we have $\mathcal{E}'_s(\R^n)\subset \mathcal{D}'_s(\R^n)$ and $\mathcal{E}'_{(s)}(\R^n)\subset \mathcal{D}'_{(s)}(\R^n)$.
We also have $\mathcal{D}'(\R^n)\subset \mathcal{D}'_s(\R^n)\subset \mathcal{D}'_{(s)}(\R^n)$.
The Fourier transform of the functionals of $\mathcal{E}'_s(\R^n)$ and
$\mathcal{E}'_{(s)}(\R^n)$ can be defined in the same way as for the distributions. Then, the following
characterisation holds (see \cite{K, LO:09, Rodino:93}):

\begin{prop}
\label{prop_ud}
A real analytic functional $v$ belongs to $\mathcal{E}'_s(\R^n)$ if and only if for any $\delta>0$ there exists $C_\delta>0$ such that
\[
|\widehat{v}(\xi)|\le C_\delta\,\esp^{\delta\lara{\xi}^{\frac{1}{s}}}
\]
for all $\xi\in\Rn$. Similarly, $v\in \mathcal{E}'_{(s)}(\R^n)$ if and only if there exist $\delta>0$ and $C>0$ such that
\[
|\widehat{v}(\xi)|\le C\,\esp^{\delta\lara{\xi}^{\frac{1}{s}}}
\]
for all $\xi\in\Rn$.
\end{prop}
We are now ready to state the ultradistributional versions of Theorem \ref{THM:case1} and
Theorem \ref{THM:case2}.

\begin{thm}\label{THM:case1u}
Let $T>0$ and $0\leq l\leq m-1$. Assume the conditions of Case 1.
Then for any $g_k\in \mathcal{E}'_{(s)}(\R^n)$ ($k=1,\ldots,m$),
the Cauchy problem \eqref{CP} has a unique global solution
$u\in C^m([0,T];\mathcal{D}'_{(s)}(\R^n))$, provided that
$$
1\leq s\leq 1+\min\left\{\alpha,\frac{m-l}{l}\right\}.
$$
\end{thm}
The situation in Case 2 is as follows:

\begin{thm}\label{THM:case2u}
Let $T>0$, $2\leq r\leq m-1$ and $0\leq l\leq m-1$. Assume the conditions of Case 2.
Then for any $g_k\in \mathcal{E}'_{(s)}(\R^n)$ ($k=1,\ldots,m$),
the Cauchy problem \eqref{CP} has a unique global solution
$u\in C^m([0,T]; \mathcal{D}'_{(s)}(\R^n))$, provided that
$$
1\leq s\leq 1+\min\left\{\alpha,\frac{\beta}{r-\beta}\right\}.
$$
\end{thm}
It is interesting to note
the non-strict inequalities for $s$ in Theorems
\ref{THM:case1u} and \ref{THM:case2u} as opposed to
strict inequalities for $s$ in Theorems \ref{THM:case1}
and \ref{THM:case2}, see also Remark
\ref{REM:local}.

Finally, we make a remark about the strictly hyperbolic case with low
regularity coefficients.
\begin{rem}\label{REM:r1statemement-u}
Under the conditions of Case 3, the conclusion of
Theorem  \ref{THM:case2u} holds provided that $$1\le s<1+\frac{\beta}{1-\beta}.$$
See Remark \ref{REM:r1} for the argument.
\end{rem}

\section{Reduction to first order system and preliminary analysis}

We now perform a reduction to a first order system as in \cite{Taylor:81}. Let  $\lara{D_x}$ be
the pseudo-differential operator with symbol $\lara{\xi}$. The transformation
\[
u_k=D_t^{k-1}\lara{D_x}^{m-k}u,
\]
with $k=1,...,m$, makes the Cauchy problem \eqref{CP} equivalent to the following system
\beq
\label{syst_Taylor}
D_t\left(
                             \begin{array}{c}
                               u_1 \\
                               \cdot \\
                               \cdot\\
                               u_m \\
                             \end{array}
                           \right)
= \left(
    \begin{array}{ccccc}
      0 & \lara{D_x} & 0 & \dots & 0\\
      0 & 0 & \lara{D_x} & \dots & 0 \\
      \dots & \dots & \dots & \dots & \lara{D_x} \\
      b_1 & b_2 & \dots & \dots & b_m \\
    \end{array}
  \right)
  \left(\begin{array}{c}
                               u_1 \\
                               \cdot \\
                               \cdot\\
                               u_m \\
                             \end{array}
                           \right)
                           +
\left(\begin{array}{c}
                               0 \\
                               0 \\
                               \cdot\\
                               f \\
                             \end{array}
                           \right),
\eeq
where
\[
b_j=A_{m-j+1}(t,D_x)\lara{D_x}^{j-m},
\]
with initial condition
\beq
\label{ic_Taylor}
u_k|_{t=0}=\lara{D_x}^{m-k}g_k,\qquad k=1,...,m.
\eeq
The matrix in \eqref{syst_Taylor} can be written as $A+B$ with
\[
A=\left(
    \begin{array}{ccccc}
      0 & \lara{D_x} & 0 & \dots & 0\\
      0 & 0 & \lara{D_x} & \dots & 0 \\
      \dots & \dots & \dots & \dots & \lara{D_x} \\
      b_{(1)} & b_{(2)} & \dots & \dots & b_{(m)} \\
    \end{array}
  \right),
\]
where $b_{(j)}=A_{(m-j+1)}(t,D_x)\lara{D_x}^{j-m}$ and
\[
B=\left(
    \begin{array}{ccccc}
      0 & 0 & 0 & \dots & 0\\
      0 & 0 & 0& \dots & 0 \\
      \dots & \dots & \dots & \dots & 0 \\
      b_1-b_{(1)} & b_2-b_{(2)} & \dots & \dots & b_m-b_{(m)} \\
    \end{array}
  \right).
\]
It is clear that the eigenvalues of the symbol matrix $A(t,\xi)$ are the roots $\tau_j(t,\xi)$, $j=1,...,m$.
By Fourier transforming both sides of \eqref{syst_Taylor} we obtain the system
\beq
\label{system_new}
\begin{split}
D_t V&=A(t,\xi)V+B(t,\xi)V+\widehat{F}(t,\xi),\\
V|_{t=0}(\xi)&=V_0(\xi),
\end{split}
\eeq
where $V$ is the $m$-column with entries $v_k=\widehat{u}_k$, $V_0$ is the $m$-column with entries
$v_{0,k}=\lara{\xi}^{m-k}\widehat{g}_k$ and
\begin{multline*}
A(t,\xi)=\left(
    \begin{array}{ccccc}
      0 & \lara{\xi} & 0 & \dots & 0\\
      0 & 0 & \lara{\xi} & \dots & 0 \\
      \dots & \dots & \dots & \dots & \lara{\xi} \\
      b_{(1)}(t,\xi) & b_{(2)}(t,\xi) & \dots & \dots & b_{(m)}(t,\xi) \\
    \end{array}
  \right),\\[0.3cm]
  b_{(j)}(t,\xi)=A_{(m-j+1)}(t,\xi)\lara{\xi}^{j-m},
\end{multline*}
\begin{multline*}
B(t,\xi)=\left(
    \begin{array}{ccccc}
      0 & 0 & 0 & \dots & 0\\
      0 & 0 & 0& \dots & 0 \\
      \dots & \dots & \dots & \dots & 0 \\
      (b_1-b_{(1)})(t,\xi) & \dots & \dots & \dots & (b_m-b_{(m)})(t,\xi) \\
    \end{array}
  \right),\\[0.3cm]
(b_j-b_{(j)})(t,\xi)=(A_{m-j+1}-A_{(m-j+1)})(t,\xi)\lara{\xi}^{j-m},
\end{multline*}
\[
\widehat{F}(t,\xi)=\left(\begin{array}{c}
                               0 \\
                               0 \\
                               \vdots\\
                               \widehat{f}(t,\cdot)(\xi) \\
                             \end{array}
                           \right).
\]
From now on we will concentrate on the system \eqref{system_new}.
We collect some preliminary results which will be crucial in the next section. Detailed proofs can be obtained by easily adapting the Lemmas 1, 2, 4 and 5 in \cite[Section 2]{ColKi:02} to our situation.
\begin{prop}
\label{prop_prelim}
Let $\lambda_i\in\R$, $i=1,...,m$, be distinct and let
\beq
\label{def_H}
H=\left(
    \begin{array}{ccccc}
      1 & 1 & 1 & \dots & 1\\
      \lambda_1\lara{\xi}^{-1} & \lambda_2\lara{\xi}^{-1} & \lambda_3\lara{\xi}^{-1} & \dots & \lambda_m\lara{\xi}^{-1} \\
      \lambda^2_1\lara{\xi}^{-2} & \lambda^2_2\lara{\xi}^{-2} & \lambda^2_3\lara{\xi}^{-2} & \dots & \lambda^2_m\lara{\xi}^{-2} \\
      \dots & \dots & \dots & \dots & \dots\\
      \lambda^{m-1}_1\lara{\xi}^{-m+1} & \lambda^{m-1}_2\lara{\xi}^{-m+1} & \lambda^{m-1}_3\lara{\xi}^{-m+1} & \dots & \lambda^{m-1}_m\lara{\xi}^{-m+1}\\
    \end{array}
  \right).
\eeq
Then we have the following properties:
\begin{itemize}
\item[(i)] $\det H=\lara{\xi}^{-\frac{(m-1)m}{2}}\prod_{1\le j<i\le m}(\lambda_i-\lambda_j)$ and
\[
\det (A(t,\xi)-\tau I)=(-1)^m(\tau^m-\sum_{j=0}^{m-1}A_{(m-j)}(t,\xi)\tau^{j});
\]
\item[(ii)] the matrix $H^{-1}$ has entries $h_{pq}$ as follows:
\[
h_{pq}= (-1)^{q-1}\lara{\xi}^{q-1}\sum_{S^{(m)}_p(m-q)}\lambda_{i_1}\dots\lambda_{i_{m-q}}
\biggl(\prod_{i=1,i\neq p}^m (\lambda_i-\lambda_p)\biggr)^{-1},
\]
for $1\le q\le m-1$,
and
\[
h_{pq}=(-1)^{m-1}\lara{\xi}^{m-1}\biggl(\prod_{i=1,i\neq p}^m (\lambda_i-\lambda_p)\biggr)^{-1},
\]
for $q=m$,
where $$S^{(a)}_b(c)=\{(i_1,...,i_c)\in \N^c; 1\le i_1<\cdots <i_c\le a, i_k\neq b, 1\le k\le c\}.$$
\item[(iii)] the matrix $H^{-1}A(t,\xi)H$ has entries
\[
c_{pq}=(\tau_q-\lambda_q)\frac{\prod_{i=1, i\neq q}^m  (\tau_i-\lambda_q)}{\prod_{i=1, i\neq p}^m (\lambda_i-\lambda_p)}
\]
when $p\neq q$.
\item[(iv)]  the matrix $H^{-1}B(t,\xi)H$ has entries
\[
d_{pq}=(-1)^{m-1}\biggl(\prod_{i=1,i\neq p}^m (\lambda_i-\lambda_p)\biggr)^{-1}g(\lambda_q),
\]
where $g(\tau)=\sum_{j=0}^{m-1}(A_{m-j}-A_{(m-j)})(t,\xi)\tau^{j}$.
\item[(v)] Assume that $\lambda_j\in C^1(\R_t)$, $j=1,...,m$. The matrix $H^{-1}\frac{d}{dt} H$ has entries
\[
e_{pq}=\begin{cases}
-\lambda'_p(t)\sum_{i=1, i\neq p}^m\frac{1}{\lambda_i(t)-\lambda_p(t)},& p=q,\\[0.3cm]
-\lambda'_q(t)\frac{\prod_{i=1, i\neq p,q}^{m}(\lambda_i(t)-\lambda_q(t))}{\prod_{i=1, i\neq p}^{m}(\lambda_i(t)-\lambda_p(t))},& p\neq q.
\end{cases}
\]
\end{itemize}
\end{prop}
\begin{pf}
We only prove assertions (iii) and (iv) and (v).

(iii) Let $w(\tau)=\sum_{j=0}^{m-1}A_{(m-j)}(t,\xi)\tau^{j}$. Arguing as in the proof of Lemma 5 in \cite{ColKi:02} we have that
\begin{multline*}
(c_{pq})_{1\le p,q\le m}= \\[0.3cm] H^{-1}\left(
    \begin{array}{ccccc}
      \lambda_1 & \lambda_2 & \lambda_3 & \dots & \lambda_m\\
      \lambda_1^2\lara{\xi}^{-1} & \lambda^2_2\lara{\xi}^{-1} & \lambda^2_3\lara{\xi}^{-1} & \dots & \lambda^2_m\lara{\xi}^{-1} \\
      \lambda^3_1\lara{\xi}^{-2} & \lambda^3_2\lara{\xi}^{-2} & \lambda^3_3\lara{\xi}^{-2} & \dots & \lambda^3_m\lara{\xi}^{-2} \\
      \dots & \dots & \dots & \dots & \dots \\
      w(\lambda_1)\lara{\xi}^{-m+1} & w(\lambda_2)\lara{\xi}^{-m+1} & w(\lambda_3)\lara{\xi}^{-m+1} & \dots & w(\lambda_m)\lara{\xi}^{-m+1}\\
    \end{array}
  \right).
\end{multline*}
Assertion (ii) yields
\begin{multline*}
c_{pq}=\sum_{r=1}^{m-1}h_{pr}\lambda_q^r\lara{\xi}^{-r+1}+h_{pm}\lara{\xi}^{-m+1}f(\lambda_q)\\
=\sum_{r=1}^{m-1}(-1)^{r-1}\sum_{S^{(m)}_p(m-r)}\lambda_{i_1}\dots\lambda_{i_{m-r}}\biggl(\prod_{i=1,i\neq p}^m (\lambda_i-\lambda_p)\biggr)^{-1}\lambda_q^r\\
+(-1)^{m-1}\biggl(\prod_{i=1,i\neq p}^m (\lambda_i-\lambda_p)\biggr)^{-1}f(\lambda_q),
\end{multline*}
which coincides with formula (25) in \cite{ColKi:02}. The proof continues as in \cite[Lemma 5]{ColKi:02}.

(iv) Let $g(\tau)=\sum_{j=0}^{m-1}(A_{m-j}-A_{(m-j)})(t,\xi)\tau^{j}$. The matrix $H^{-1}B(t,\xi)H$ can be  written as
\begin{multline*}
(d_{pq})_{1\le p,q\le m}= \\[0.3cm]
H^{-1}\left(
    \begin{array}{ccccc}
      0 & 0 & 0 & \dots & 0\\
      0 & 0 & 0 & \dots & 0\\
      0 & 0 & 0 & \dots & 0 \\
      \dots & \dots & \dots & \dots & \dots \\
      g(\lambda_1)\lara{\xi}^{-m+1} & g(\lambda_2)\lara{\xi}^{-m+1} & g(\lambda_3)\lara{\xi}^{-m+1} & \dots & g(\lambda_m)\lara{\xi}^{-m+1}\\
    \end{array}
  \right).
\end{multline*}
From (ii) we conclude that
\begin{multline*}
d_{pq}=(-1)^{m-1}\lara{\xi}^{m-1}\biggl(\prod_{i=1,i\neq p}^m (\lambda_i-\lambda_p)\biggr)^{-1}\lara{\xi}^{-m+1}g(\lambda_q)\\
= (-1)^{m-1}\biggl(\prod_{i=1,i\neq p}^m (\lambda_i-\lambda_p)\biggr)^{-1}g(\lambda_q).
\end{multline*}

(v) From the definition of $H$ we have that $H^{-1}\frac{d}{dt} H$ is the matrix
\[
H^{-1}\left(
    \begin{array}{ccccc}
      0 & 0 & 0 & \dots & 0\\
      \lambda'_1\lara{\xi}^{-1} & \lambda'_2\lara{\xi}^{-1}& \lambda'_3\lara{\xi}^{-1}& \dots & \lambda'_m\lara{\xi}^{-1}\\
      (\lambda^2_1)'\lara{\xi}^{-2}& (\lambda^2_2)'\lara{\xi}^{-2} & (\lambda^2_3)'\lara{\xi}^{-2} & \dots & (\lambda^2_m)'\lara{\xi}^{-2} \\
      \dots & \dots & \dots & \dots & \dots \\
      (\lambda^{m-1}_1)'\lara{\xi}^{-m+1}& (\lambda^{m-1}_2)'\lara{\xi}^{-m+1} & (\lambda^{m-1}_3)'\lara{\xi}^{-m+1} & \dots & (\lambda^{m-1}_m)'\lara{\xi}^{-m+1}\\
    \end{array}
  \right).
\]
Hence, making use of the second assertion of this proposition we obtain
\begin{multline*}
e_{pq}=\sum_{r=2}^{m-1}h_{pr}(r-1)\lambda^{r-2}_q\lambda'_q\lara{\xi}^{-r+1}+h_{pm}(m-1)\lambda^{m-2}_q\lambda'_q
\lara{\xi}^{-m+1}\\
=\sum_{r=2}^{m-1}(-1)^{r-1}\lara{\xi}^{r-1}\sum_{S^{(m)}_p(m-r)}\lambda_{i_1}\dots\lambda_{i_{m-r}}\biggl(\prod_{i=1,i\neq p}^m (\lambda_i-\lambda_p)\biggr)^{-1}\lambda^{r-2}_q\lambda'_q\lara{\xi}^{-r+1}\\
+(-1)^{m-1}\lara{\xi}^{m-1}\biggl(\prod_{i=1,i\neq p}^m (\lambda_i-\lambda_p)\biggr)^{-1}(m-1)\lambda^{m-2}_q\lambda'_q
\lara{\xi}^{-m+1}\\
=\sum_{r=2}^{m-1}(-1)^{r-1}\sum_{S^{(m)}_p(m-r)}\lambda_{i_1}\dots\lambda_{i_{m-r}}\biggl(\prod_{i=1,i\neq p}^m (\lambda_i-\lambda_p)\biggr)^{-1}\lambda^{r-2}_q\lambda'_q\\
+(-1)^{m-1}\biggl(\prod_{i=1,i\neq p}^m (\lambda_i-\lambda_p)\biggr)^{-1}(m-1)\lambda^{m-2}_q\lambda'_q.
\end{multline*}
This is the expression for $b_{pq}$ in the proof of Lemma 4 in \cite{ColKi:02}. The proof continues as in \cite[Lemma 4]{ColKi:02}.
\end{pf}
We now proceed to analyse the roots $\tau_j$. We perform the natural
regularisation and separation process, but it will be different under the assumptions of Case 1 or of Case 2.
To simplify the notation, although the functions below will depend on
$\eps$, for brevity we will write
$\lambda_j(t,\xi)$ for $\lambda_j(\eps,t,\xi)$.
\begin{prop}
\label{prop_roots}
Let $\varphi\in C^\infty_{c}(\R)$, $\varphi\ge 0$ with $\int_\R\varphi(x)\, dx=1$. Under the assumptions of Case 1, let
\beq
\label{def_lambdaj}
\lambda_j(t,\xi)=(\tau_j(\cdot,\xi)\ast\varphi_\eps)(t)+j\eps^\alpha\lara{\xi},
\eeq
for $j=1,...,m$ and $\varphi_\eps(s)=\eps^{-1}\varphi(s/\eps)$, $\eps>0$.
Then, there exists a constant $c>0$ such that
\begin{itemize}
\item[(i)] $|\partial_t\lambda_j(t,\xi)|\le c\,\eps^{\alpha-1}\lara{\xi}$,
\item[(ii)] $|\lambda_j(t,\xi)-\tau_j(t,\xi)|\le c\,\eps^{\alpha}\lara{\xi}$,
\item[(iii)] $\lambda_{j}(t,\xi)-\lambda_{i}(t,\xi)\ge \eps^\alpha\lara{\xi}$ for $j>i$,
\end{itemize}
for all $t,s\in[0,T']$ with $T'<T$ and all $\xi\in\R^n$.
\end{prop}
\begin{pf}
By definition of convolution, if $R$ is large enough, one has
\begin{multline}
\label{est_deriv}
|\partial_t\lambda_j(t,\xi)|=\eps^{-1}\int_{-R}^R\tau_j(t-\eps s)\varphi'(s)\, ds\\
=\eps^{-1}\int_{-R}^R(\tau_j(t-\eps s,\xi)-\tau_j(t,\xi))\varphi'(s)\, ds+\eps^{-1}\int_{-R}^R \tau_j(t,\xi)\varphi'(s)\, ds,
\end{multline}
and, therefore, by \eqref{hyp_roots} we obtain
$|\partial_t\lambda_j(t,\xi)|\le c\eps^{\alpha-1}\lara{\xi}$ for all $t,s\in[0,T']$ and $\xi\in\R^n$.
 The second and third assertions follow immediately from the definition of $\lambda_j$, where we note that in view of
\eqref{hyp_roots_2} and the fact that
$\varphi\geq 0$ it is enough to observe (iii) for $j-i=1$.
\end{pf}
\begin{prop}
\label{prop_roots_2}
Let $\varphi\in C^\infty_{c}(\R)$, $\varphi\ge 0$ with $\int_\R\varphi(x)\, dx=1$. Under the assumptions of Case 2, let
\beq
\label{def_lambdaj_2}
\begin{split}
\lambda_j(t,\xi)&=(\tau_j(\cdot,\xi)\ast\varphi_\eps)(t)+j\eps^\alpha\lara{\xi},\quad 1\le j\le r,\\
\lambda_j(t,\xi)&=(\tau_j(\cdot,\xi)\ast\varphi_\delta)(t),\quad\qquad r+1\le j\le m,
\end{split}
\eeq
for $0<\delta,\eps<1$. Then, there exist constants $c>0$, $c_0>0$ such that
\begin{itemize}
\item[(i)] $|\partial_t\lambda_j(t,\xi)|\le c\,\eps^{\alpha-1}\lara{\xi}$ for $j=1,...,r$,
\item[(ii)] $|\lambda_j(t,\xi)-\tau_j(t,\xi)|\le c\,\eps^{\alpha}\lara{\xi}$ for $j=1,...,r$,
 \item[(iii)] $\lambda_{j+1}(t,\xi)-\lambda_{j}(t,\xi)\ge \eps^\alpha\lara{\xi}$ for $j=1,...,r-1$,
\item[(iv)] $|\partial_t\lambda_j(t,\xi)|\le c\,\delta^{\beta-1}\lara{\xi}$ for $j=r+1,...,m$,
\item[(v)] $|\lambda_{j}(t,\xi)-\tau_j(t,\xi)|\le c\,\delta^{\beta}\lara{\xi}$ for $j=r+1,...,m$,
\item[(vi)] $\lambda_{j+1}(t,\xi)-\lambda_j(t,\xi)\ge c_0\lara{\xi}$ for $j=r,...,m-1$, for
$\eps=\lara{\xi}^{-\gamma}$ with $\gamma\in(0,1)$, $\delta=\lara{\xi}^{-1}$ and $|\xi|$ large enough,
\item[(vii)] $\lambda_{j}(t,\xi)-\lambda_i(t,\xi)\ge c_0\lara{\xi}$
for $j=r+1,...,m$, $i=1,...,r$, $\eps=\lara{\xi}^{-\gamma}$ with $\gamma\in(0,1)$, $\delta=\lara{\xi}^{-1}$ and $|\xi|$ large enough,
\end{itemize}
hold for all $t,s\in[0,T']$ with $T'<T$.
\end{prop}
\begin{pf}
The first three assertions are clear from Proposition \ref{prop_roots} and
\eqref{hyp_roots_3}. Assertion (iv) can be proven as in \eqref{est_deriv}.
Assertion (v) follows immediately from the $C^\beta$-property of the roots $\tau_j$ when $j=r+1,...,m$. We finally consider the difference $\lambda_{j+1}(t,\xi)-\lambda_j(t,\xi)$. If $j=r+1,...,m-1$ then from the bound from below \eqref{strict_hyp} we obtain the estimate
\[
\lambda_{j+1}(t,\xi)-\lambda_j(t,\xi)\ge c_0\lara{\xi}
\]
valid for $t\in[0,T']$ and $|\xi|$ large enough. It remains to consider $\lambda_{j+1}(t,\xi)-\lambda_j(t,\xi)$ when $j=r$. Making use of the definition in \eqref{def_lambdaj_2} we can write
\begin{multline*}
\lambda_{r+1}(t,\xi)-\lambda_r(t,\xi)=\int_\R \tau_{r+1}(t-\delta s,\xi)\varphi(s)\, ds -\int_\R \tau_{r}(t-\eps s,\xi)\varphi(s)\, ds-r\eps^\alpha\lara{\xi}\\
= \int_\R (\tau_{r+1}(t-\delta s,\xi)-\tau_{r+1}(t-\eps s,\xi))\varphi(s)\, ds +\\
 +\int_\R (\tau_{r+1}(t-\eps s,\xi)-\tau_r(t-\eps s,\xi))\varphi(s)\, ds -r\eps^\alpha\lara{\xi}.
\end{multline*}
Hence, combining \eqref{strict_hyp} with \eqref{hyp_roots} we get
\[
\lambda_{r+1}(t,\xi)-\lambda_r(t,\xi)\ge c_0|\xi|-c|\eps-\delta|^\beta|\xi|-r\eps^\alpha\lara{\xi}\ge c_0|\xi|-c|\eps-\delta|^\beta|\xi|-r\eps^\alpha\sqrt{2}|\xi|,
\]
for $|\xi|\ge 1$. It follows that for
\beq
\label{radius_1}
|\eps-\delta|^\beta\le \frac{c_0}{4c} \Leftrightarrow |\eps-\delta|\le \big({\frac{c_0}{4c}}\big)^{\frac{1}{\beta}} \Leftrightarrow \lara{\xi}^{-\gamma}(1-\lara{\xi}^{-1+\gamma})\le \big({\frac{c_0}{4c}}\big)^{\frac{1}{\beta}}
\eeq
and
\beq
\label{radius_2}
\eps^\alpha\le \frac{c_0}{4\sqrt{2}r} \Leftrightarrow \eps\le \big({\frac{c_0}{4\sqrt{2}r}}\big)^{\frac{1}{\alpha}} \Leftrightarrow \lara{\xi}^{-\gamma}\le \big({\frac{c_0}{4\sqrt{2}r}}\big)^{\frac{1}{\alpha}}
\eeq
one has
\[
\lambda_{r+1}(t,\xi)-\lambda_r(t,\xi)\ge c'_0\lara{\xi}.
\]
%This bound from below is true for $|\xi|\ge R$, with $R$ determined by \eqref{radius_1} and \eqref{radius_2}.
%Since $\varphi_\eps$ is an approximation of the Dirac delta distribution, it follows that for
%$\lara{\xi}\to\infty$ (or for $\eps\to 0$), we get that
%$\lambda_{r+1}(t,\xi)-\lambda_r(t,\xi)$
%approximates the positive quantity $\tau_{r+1}(t,\xi)-\tau_r(t,\xi)>0.$ Consequently,
%we can remove the modulus in (vi).

Assertion (vii) follows from (vi).
\end{pf}
In the sequel, with abuse of notation, we will still denote the smaller $T'$ in Propositions
\ref{prop_roots} and \ref{prop_roots_2} by $T$.
\begin{prop}
\label{prop_roots_3}
The property \eqref{hyp_coincide} holds for the $\lambda_j$'s as well, i.e.,
\beq
\label{hyp_coincide_2}
|\lambda_i(t,\xi)-\lambda_j(t,\xi)|\le c|\lambda_k(t,\xi)-\lambda_{k-1}(t,\xi)|
\eeq
for all $1\le i,j,k\le r$, for all $t\in[0,T]$ and $\xi\in\R^n$.
\end{prop}
\begin{pf}
Assume that $i>j$. Hence
\[
|\lambda_i(t,\xi)-\lambda_j(t,\xi)|=(\tau_i(\cdot,\xi)-\tau_j(\cdot,\xi))\ast\varphi_\eps(t)+
(i-j)\eps^\alpha\lara{\xi}
\]
and
\[
|\lambda_k(t,\xi)-\lambda_{k-1}(t,\xi)|=(\tau_{k}(\cdot,\xi)-\tau_{k-1}(\cdot,\xi))\ast\varphi_\eps(t)+
\eps^\alpha\lara{\xi}.
\]
From \eqref{hyp_coincide} and the fact that $\varphi\geq 0$ we get that
\begin{multline*}
|\lambda_i(t,\xi)-\lambda_j(t,\xi)|\le c(\tau_{k}(\cdot,\xi)-\tau_{k-1}
(\cdot,\xi))\ast\varphi_\eps(t)+(i-j)\eps^\alpha\lara{\xi}\\[0.2cm]
\le c'|\lambda_k(t,\xi)-\lambda_{k-1}(t,\xi)|
\end{multline*}
holds for all $t\in[0,T]$ and $\xi\in\R^n$.
\end{pf}

%\begin{remark}
%Note that similarly to assertion $(vi)$ one has that
%\[
%|\lambda_{j}(t,\xi)-\lambda_i(t,\xi)|\ge c_0\lara{\xi}
%\]
%is valid for $j=r+1,...,m$, $i=1,...,r$, $t\in[0,T']$ and $|\xi|$ large enough.
%\end{remark}

\section{Proof in Case 1: Theorem \ref{THM:case1} and Theorem \ref{THM:case1u}}

We first prove Theorem  \ref{THM:case1}.
It is well-known that the problem \eqref{syst_Taylor}-\eqref{ic_Taylor} is well-posed when $s=1$,
see e.g. \cite{Jan:84, Kaj:86}. Hence, we may assume $s>1$.
In the case of Theorem  \ref{THM:case1} we can also assume that the initial data have compact support. Since weakly hyperbolic equations have the finite speed of propagation property it follows that the solution $u$ is compactly supported in $x$ as well. This observation allows us to proceed with the reduction to a first order system of Section 3.

Let $H(t,\xi)$ be the matrix \eqref{def_H} with entries $\lambda_j(t,\xi)$ as in \eqref{def_lambdaj}. Observe that the approximated roots $\lambda_j$ are distinct for all $\eps>0$. We look for a solution $V$ of the Cauchy problem \eqref{system_new} in the form
\beq
\label{VW}
V(t,\xi)=\esp^{-\rho(t)\lara{\xi}^{\frac{1}{s}}}(\det H)^{-1}HW,
\eeq
where $\rho\in C^1[0,T]$ will be determined in the sequel. By substitution in \eqref{system_new} we obtain
\begin{multline*}
\esp^{-\rho(t)\lara{\xi}^{\frac{1}{s}}}(\det H)^{-1}HD_tW+\esp^{-\rho(t)\lara{\xi}^{\frac{1}{s}}}\irm\rho'(t)\lara{\xi}^{\frac{1}{s}}(\det H)^{-1}HW+\\
+\irm\esp^{-\rho(t)\lara{\xi}^{\frac{1}{s}}}\frac{\partial_t\det H}{(\det H)^2}HW +\esp^{-\rho(t)\lara{\xi}^{\frac{1}{s}}}(\det H)^{-1}(D_tH)W\\
=\esp^{-\rho(t)\lara{\xi}^{\frac{1}{s}}}(\det H)^{-1}(A+B)HW + \widehat{F}.
\end{multline*}
Multiplying both sides of the previous equation by $\esp^{\rho(t)\lara{\xi}^{\frac{1}{s}}}(\det H)H^{-1}$ we get
\begin{multline*}
D_tW+\irm\rho'(t)\lara{\xi}^{\frac{1}{s}}W+\irm\frac{\partial_t\det H}{\det H}W + H^{-1}(D_t H)W= H^{-1}(A+B)HW+\\
+\esp^{\rho(t)\lara{\xi}^{\frac{1}{s}}}(\det H) H^{-1}\widehat{F}.
\end{multline*}
Hence,
\begin{multline}
\label{energy}
\partial_t |W(t,\xi)|^2=2{\rm Re} (\partial_t W(t,\xi),W(t,\xi))\\
=2\rho'(t)\lara{\xi}^{\frac{1}{s}}|W(t,\xi)|^2+2\frac{\partial_t\det H}{\det H}|W(t,\xi)|^2-2
{\rm Re}(H^{-1}\partial_t HW,W)\\
-2{\rm Im} (H^{-1}AHW,W)-2{\rm Im} (H^{-1}BHW,W) \\
-2{\rm Im} (\esp^{\rho(t)\lara{\xi}^{\frac{1}{s}}}(\det H) H^{-1}\widehat{F},W).
\end{multline}
We proceed by estimating
\begin{enumerate}
\item $\frac{\partial_t\det H}{\det H}$,
\item $\Vert H^{-1}\partial_t H\Vert$,
\item $\Vert H^{-1}AH-(H^{-1}AH)^\ast\Vert$,
\item $\Vert H^{-1}BH-(H^{-1}BH)^\ast\Vert$.
\end{enumerate}
\subsection{Estimate of the first term}
Proposition \ref{prop_prelim}(i) combined with Proposition \ref{prop_roots} yields the following estimate
\begin{multline}
\label{est_1}
\biggl|\frac{\partial_t\det H(t,\xi)}{\det H(t,\xi)}\biggr|=\biggl|
\frac{\lara{\xi}^{-\frac{(m-1)m}{2}}\partial_t\prod_{1\le j<i\le m}
(\lambda_i(t,\xi)-\lambda_j(t,\xi))}{\lara{\xi}^{-\frac{(m-1)m}{2}}
\prod_{1\le j<i\le m}(\lambda_i(t,\xi)-\lambda_j(t,\xi))}\biggr| \\
\le \sum_{1\le j<i\le  m}\frac{|\partial_t\lambda_i(t,\xi)-
\partial_t\lambda_j(t,\xi)|}{|\lambda_i(t,\xi)-\lambda_j(t,\xi)|}
\le\frac{c_1\eps^{\alpha-1}\lara{\xi}}{\eps^\alpha\lara{\xi}}=c_1\,\eps^{-1},
\end{multline}
valid for all $t\in[0,T]$ and $\xi\in\R^n$.
\subsection{Estimate of the second term}
\label{second_1}
From Proposition \ref{prop_prelim}(v) the entries of the matrix $H^{-1}(t,\xi)\partial_t H(t,\xi)$ can be written as
\[
e_{pq}(t,\xi)=\begin{cases}
-\partial_t\lambda_p(t,\xi)\sum_{i=1, i\neq p}^m\frac{1}{\lambda_i(t,\xi)-\lambda_p(t,\xi)},& p=q,\\[0.3cm]
-\partial_t\lambda_q(t,\xi)\frac{\prod_{i=1, i\neq p,q}^{m}(\lambda_i(t,\xi)-\lambda_q(t,\xi))}{\prod_{i=1, i\neq p}^{m}(\lambda_i(t,\xi)-\lambda_p(t,\xi))},& p\neq q.
\end{cases}
\]
From Proposition \ref{prop_roots} we clearly have that
\[
|e_{pp}(t,\xi)|\le c\frac{\eps^{\alpha-1}\lara{\xi}}{\eps^\alpha\lara{\xi}}=c\eps^{-1}.
\]
To estimate $e_{pq}$ when $q\neq p$ we write
\[
\partial_t\lambda_q(t,\xi)\frac{\prod_{i=1, i\neq p,q}^{m}(\lambda_i(t,\xi)-\lambda_q(t,\xi))}{\prod_{i=1, i\neq p}^{m}(\lambda_i(t,\xi)-\lambda_p(t,\xi))}
\]
as
\[
\partial_t\lambda_q(t,\xi)\frac{\prod_{i=1, i\neq p,q}^{m}(\lambda_i(t,\xi)-\lambda_q(t,\xi))}{\prod_{i=1, i\neq p,q}^{m}(\lambda_i(t,\xi)-\lambda_p(t,\xi))(\lambda_q(t,\xi)-\lambda_p(t,\xi))}.
\]
Since $$|\lambda_i(t,\xi)-\lambda_q(t,\xi)|\le |\lambda_i(t,\xi)-\lambda_p(t,\xi)|+
|\lambda_p(t,\xi)-\lambda_q(t,\xi)|$$ arguing as in (40) in \cite{ColKi:02} and making use of the estimate \eqref{hyp_coincide_2} we obtain that
\[
|e_{pq}(t,\xi)|\le c\frac{\eps^{\alpha-1}\lara{\xi}}{\eps^\alpha\lara{\xi}}=c\eps^{-1}.
\]
Hence, $\Vert H^{-1}\partial_t H\Vert\le c_2\eps^{-1}$.

\subsection{Estimate of the third term}
From Proposition \ref{prop_prelim}(iii) the matrix $H^{-1}AH$ has entries
\[
c_{pq}(t,\xi)=(\tau_q(t,\xi)-\lambda_q(t,\xi))\frac{\prod_{i=1, i\neq q}^m  (\tau_i(t,\xi)-\lambda_q(t,\xi))}{\prod_{i=1, i\neq p}^m (\lambda_i(t,\xi)-\lambda_p(t,\xi))}
\]
when $p\neq q$. As in formula (46) in \cite{ColKi:02} we have
\begin{multline*}
|\tau_q(t,\xi)-\lambda_q(t,\xi)|\frac{\prod_{i=1, i\neq q}^m  |\tau_i(t,\xi)-\lambda_q(t,\xi)|}{\prod_{i=1, i\neq p}^m |\lambda_i(t,\xi)-\lambda_p(t,\xi)|}\\
\le |\tau_q(t,\xi)-\lambda_q(t,\xi)|\frac{\prod_{i=1, i\neq q}^m  |\tau_q(t,\xi)-\lambda_q(t,\xi)|+|\tau_i(t,\xi)-\tau_q(t,\xi)|}{\prod_{i=1, i\neq p}^m |\lambda_i(t,\xi)-\lambda_p(t,\xi)|}\\
=\sum_{k=1}^{m-1}|\tau_q(t,\xi)-\lambda_q(t,\xi)|^k\sum_{S^{(m)}_q(m-k)}\frac{|\tau_{i_1}(t,\xi)-\tau_q(t,\xi)|\cdots
|\tau_{i_{m-k}}(t,\xi)-\tau_q(t,\xi)|}{\prod_{i=1, i\neq p}^m |\lambda_i(t,\xi)-\lambda_p(t,\xi)|}\\
+\frac{|\tau_q(t,\xi)-\lambda_q(t,\xi)|^m}{\prod_{i=1, i\neq p}^m |\lambda_i(t,\xi)-\lambda_p(t,\xi)|}.
\end{multline*}
Proposition \ref{prop_roots} combined with
\[
|\tau_{i_k}(t,\xi)-\tau_q(t,\xi)|\le |\tau_{i_k}(t,\xi)-\lambda_{i_k}(t,\xi)|+|\lambda_{i_k}(t,\xi)-\lambda_q(t,\xi)|+|\lambda_q(t,\xi)-\tau_q(t,\xi)|,
\]
the property \eqref{hyp_coincide} and the fact that $|\tau_i(t,\xi)-\tau_j(t,\xi)|/|\lambda_i(t,\xi)-\lambda_j(t,\xi)|$ is bounded when $i\neq j$, yields the estimate
\[
|c_{pq}(t,\xi)|\le c\sum_{k=1}^{m-1}\eps^{\alpha k}\lara{\xi}^k\sum_{S^{(m)}_q(m-k)}\frac{\lara{\xi}^{m-k-m+1}}{\eps^{\alpha(m-1)-\alpha(m-k)}}+c\,\frac{\eps^{\alpha m}\lara{\xi}^m}{\eps^{\alpha(m-1)}\lara{\xi}^{m-1}}\le c\, \eps^\alpha\lara{\xi}.
\]
This implies $\Vert H^{-1}AH-(H^{-1}AH)^\ast\Vert\le c_3\eps^\alpha\lara{\xi}$.

\subsection{Estimate of the fourth term}
From Proposition \ref{prop_prelim}(iv) we have that $H^{-1}BH$ has entries
\[
d_{pq}(t,\xi)=(-1)^{m-1}\biggl(\prod_{i=1,i\neq p}^m (\lambda_i(t,\xi)-\lambda_p(t,\xi))\biggr)^{-1}g(\lambda_q(t,\xi)),
\]
where
$$g(\tau)=\sum_{j=0}^{m-1}(A_{m-j}-A_{(m-j)})(t,\xi)\tau^{j}.$$
Assume that we have lower order terms of order $l$. Then
$$|g(\lambda_q(t,\xi))|\leq C\lara{\xi}^l$$ and by Proposition \ref{prop_roots}(iii) we get
\[
|d_{pq}(t,\xi)|\le c\eps^{\alpha(1-m)}\lara{\xi}^{-m+1+l}.
\]
Hence $\Vert H^{-1}BH-(H^{-1}BH)^\ast\Vert\le c_4\eps^{\alpha(1-m)}\lara{\xi}^{l-m+1}$.

\subsection{Conclusion of the proof}
Making use of these four estimates in \eqref{energy} we get
\begin{multline}
\label{eps_energy}
\partial_t |W(t,\xi)|^2  \\
\le 2(\rho'(t)\lara{\xi}^{\frac{1}{s}}+c_1\eps^{-1}+c_2\eps^{-1}+
c_3\eps^{\alpha}\lara{\xi}+c_4\eps^{\alpha(1-m)}\lara{\xi}^{l-m+1})|W(t,\xi)|^2
\\ +  C'\esp^{(\rho(t)-\delta_1)\lara{\xi}^{\frac{1}{s}}}|W(t,\xi)|\\
\le (2\rho'(t)\lara{\xi}^{\frac{1}{s}}+C_1\eps^{-1}+C_2\eps^{\alpha}\lara{\xi}+C_3\eps^{\alpha(1-m)}\lara{\xi}^{l-m+1})
|W(t,\xi)|^2+\\
+C'\esp^{(\rho(t)-\delta_1)\lara{\xi}^{\frac{1}{s}}}|W(t,\xi)|,
\end{multline}
where $\delta_1>0$ depends on $f$, in view of Proposition \ref{prop_Fourier}.
Set $\eps=\lara{\xi}^{-\gamma}$. By substitution in \eqref{eps_energy} we arrive at comparing the terms
\[
\lara{\xi}^\gamma;\quad \lara{\xi}^{-\gamma\alpha+1};\quad \lara{\xi}^{\gamma\alpha(m-1)+l-m+1}.
\]
Choose $\gamma=\min\{\frac{1}{1+\alpha},\frac{m-l}{\alpha m}\}$. It follows that
$$\max\{\gamma, \gamma\alpha(m-1)+l-m+1\}\le -\gamma\alpha+1.$$
Then, if we take $s>0$ such that
\begin{multline}
\label{form_s}
\frac{1}{s}>-\gamma\alpha+1=-\min\biggl\{\frac{1}{1+\alpha},\frac{m-l}{\alpha m}\biggr\}\alpha+1 \\
=-\min\biggl\{\frac{\alpha}{1+\alpha},\frac{m-l}{ m}\biggr\}+1=\max\biggl\{\frac{1}{1+\alpha},\frac{l}{m}\biggr\},
\end{multline}
for a suitable decreasing function $\rho$ (for instance $\rho(t)=\rho(0)-\kappa t$ with $\kappa>0$ and $\rho(0)$ to be chosen later) we obtain
\begin{multline}\label{EQ:en-est}
\partial_t |W(t,\xi)|^2\le \big(2\rho'(t)\lara{\xi}^{\frac{1}{s}}+C\lara{\xi}^{-\gamma\alpha+1}\big)|W(t,\xi)|^2\\
+
2\esp^{\rho(t)\lara{\xi}^{\frac{1}{s}}} \det H(t,\xi) |H^{-1}(t,\xi)||\widehat{F}(t,\xi)||W(t,\xi)|\\
\le\big(2\rho'(t)\lara{\xi}^{\frac{1}{s}}+C\lara{\xi}^{-\gamma\alpha+1}\big)|W(t,\xi)|^2+
C'\esp^{(\rho(t)-\delta_1)\lara{\xi}^{\frac{1}{s}}}|W(t,\xi)|.
\end{multline}
 Note that \eqref{form_s} implies
\[
s<\min\biggl\{1+\alpha,\frac{m}{l}\biggr\}=1+\min\biggl\{\alpha, \frac{m-l}{l}\biggr\}.
\]
Assuming for the moment that $|W(t,\xi)|\ge 1$, taking $\rho(0)<\delta_1$  we get the energy estimate
\beq
\label{energy1}
\partial_t |W(t,\xi)|^2\le \big(2\rho'(t)\lara{\xi}^{\frac{1}{s}}+C\lara{\xi}^{-\gamma\alpha+1}+
C'\esp^{(\rho(0)-\delta_1)\lara{\xi}^{\frac{1}{s}}}\big)|W(t,\xi)|^2 \leq 0,
\eeq
for large enough $|\xi|$ (note that it suffices to consider only large $|\xi|$).
Consequently, \eqref{VW} and \eqref{energy1} imply the estimate
\begin{multline}
\label{last_estimate}
|V(t,\xi)|
=\esp^{-\rho(t)\lara{\xi}^{\frac{1}{s}}}\frac{1}{\det H(t,\xi)}|H(t,\xi)||W(t,\xi)|\le \\
\esp^{-\rho(t)\lara{\xi}^{\frac{1}{s}}}\frac{1}{\det H(t,\xi)}|H(t,\xi)||W(0,\xi)|=\\
\esp^{(-\rho(t)+\rho(0))\lara{\xi}^{\frac{1}{s}}}\frac{\det H(0,\xi)}{\det H(t,\xi)}|H(t,\xi)||H^{-1}(0,\xi)||V(0,\xi)|,
\end{multline}
where, for $\gamma$ as above, we have
\[
\frac{\det H(0,\xi)}{\det H(t,\xi)}|H(t,\xi)||H^{-1}(0,\xi)|\le c\,\eps^{-\alpha\frac{(m-1)m}{2}}=c\lara{\xi}^{\gamma\alpha\frac{(m-1)m}{2}}.
\]
Hence,
\beq
\label{last_estimate2}
\left\{
\begin{array}{cc}
|V(t,\xi)|\le c\,\esp^{(-\rho(t)+\rho(0))\lara{\xi}^{\frac{1}{s}}}\lara{\xi}^{\gamma\alpha\frac{(m-1)m}{2}}|V(0,\xi)|,\quad & \text{for}\, |W(t,\xi)|\ge 1,\\
|V(t,\xi)|\le c\,\esp^{-\rho(t)\lara{\xi}^{\frac{1}{s}}}\lara{\xi}^{\gamma\alpha\frac{(m-1)m}{2}},\quad &\text{for}\, |W(t,\xi)|<1,
\end{array}
\right.
\eeq
with the second line following directly from \eqref{VW}.
The estimate \eqref{last_estimate2} combined with the Fourier characterisations of Proposition \ref{prop_Fourier} yields the statement of Theorem  \ref{THM:case1} if we choose $\kappa>0$ small enough.
If $s=1+\min\biggl\{\alpha, \frac{m-l}{l}\biggr\}$, we need $\kappa$ to be large enough in
\eqref{EQ:en-est}, so that \eqref{last_estimate2} still implies the local in time
well-posedness (showing a statement in Remark \ref{REM:local}).

We note that in view of  the characterisation in Proposition \ref{prop_ud},
the estimate \eqref{last_estimate2} also yields the statement of Theorem \ref{THM:case1u}. In this case we can also allow the critical
case $s=1+\min\biggl\{\alpha, \frac{m-l}{l}\biggr\}$. Indeed, differently from
the case of Theorem \ref{THM:case1}, taking $\kappa>0$ to be large
enough, we can make sure that the estimate \eqref{energy1} holds,
while \eqref{last_estimate2} yields that $V(t,\xi)$ satisfies
the estimates of Proposition \ref{prop_ud} for any value of $T$.
Because of the presence of the function $\rho$ in \eqref{last_estimate2}
the obtained result is in the space of
Gevrey Beurling ultradistributions rather than in the space of
Gevrey Roumieu ultradistributions.

\section{Proof in Case 2:  Theorem \ref{THM:case2} and Theorem \ref{THM:case2u}}

We work on the energy estimate similar to the Case 1.
However, the different nature of the approximated roots $\lambda_j(t,\xi)$ yields different estimates for
the terms
\begin{enumerate}
\item $\frac{\partial_t\det H}{\det H}$,
\item $\Vert H^{-1}\partial_t H\Vert$,
\item $\Vert H^{-1}AH-(H^{-1}AH)^\ast\Vert$,
\item $\Vert H^{-1}BH-(H^{-1}BH)^\ast\Vert$.
\end{enumerate}

\subsection{Estimate of the first term}
Arguing as in \eqref{est_1} we have
\begin{multline*}
\biggl|\frac{\partial_t\det H(t,\xi)}{\det H(t,\xi)}\biggr|\le \sum_{1\le j<i\le  m}\frac{|\partial_t\lambda_i(t,\xi)-\partial_t\lambda_j(t,\xi)|}{|\lambda_i(t,\xi)-\lambda_j(t,\xi)|}\\
= \sum_{1\le j<i\le  r}\frac{|\partial_t\lambda_i(t,\xi)-\partial_t\lambda_j(t,\xi)|}{|\lambda_i(t,\xi)-\lambda_j(t,\xi)|}+\sum_{r+1\le j<i\le  m}\frac{|\partial_t\lambda_i(t,\xi)-\partial_t\lambda_j(t,\xi)|}{|\lambda_i(t,\xi)-\lambda_j(t,\xi)|}\\
+\sum_{\substack{1\le j<i\le m,\\ j\le r, i\ge r+1}}
\frac{|\partial_t\lambda_{i}(t,\xi)-\partial_t\lambda_j(t,\xi)|}{|\lambda_{i}(t,\xi)-\lambda_j(t,\xi)|}.
\end{multline*}
Proposition \ref{prop_roots_2} yields for $t\in[0,T]$ and $|\xi|$ large enough the following estimate:
\[
\begin{split}
\biggl|\frac{\partial_t\det H(t,\xi)}{\det H(t,\xi)}\biggr|
& \le c\frac{\eps^{\alpha-1}\lara{\xi}}{\eps^\alpha\lara{\xi}}+c'\frac{\delta^{\beta-1}\lara{\xi}}{c_0\lara{\xi}}+c''
\frac{\eps^{\alpha-1}\lara{\xi}+\delta^{\beta-1}\lara{\xi}}{c_0\lara{\xi}} \\[0.2cm]
& \le c_1\max\{\eps^{-1},\delta^{\beta-1}\}.
\end{split}
\]
We note that here we can use Proposition  \ref{prop_roots_2}(vi) since we will set
$\eps$ and $\delta$ later to be as required.

\subsection{Estimate of the second term}
The entries of the matrix $H^{-1}(t,\xi)\partial_t H(t,\xi)$ can be written as
\[
e_{pq}(t,\xi)=\begin{cases}
-\partial_t\lambda_p(t,\xi)\sum_{i=1, i\neq p}^m\frac{1}{\lambda_i(t,\xi)-\lambda_p(t,\xi)},& p=q,\\[0.3cm]
-\partial_t\lambda_q(t,\xi)\frac{\prod_{i=1, i\neq p,q}^{m}(\lambda_i(t,\xi)-\lambda_q(t,\xi))}{\prod_{i=1, i\neq p}^{m}(\lambda_i(t,\xi)-\lambda_p(t,\xi))},& p\neq q.
\end{cases}
\]
Let us start with the case $p=q$. We have
\begin{multline*}
e_{pp}(t,\xi)=-\partial_t\lambda_p(t,\xi)\sum_{i=1, i\neq p}^r\frac{1}{\lambda_i(t,\xi)-\lambda_p(t,\xi)}
\\ -\partial_t\lambda_p(t,\xi)\sum_{i=r+1, i\neq p}^m\frac{1}{\lambda_i(t,\xi)-\lambda_p(t,\xi)}.
\end{multline*}
It follows that, for $|\xi|$ large,
\[
\begin{split}
|e_{pp}(t,\xi)|&\le c\frac{\eps^{\alpha-1}\lara{\xi}}{\eps^\alpha\lara{\xi}}+
c\frac{\eps^{\alpha-1}\lara{\xi}}{c_0\lara{\xi}},\quad 1\le p\le r,\\
|e_{pp}(t,\xi)|&\le c\frac{\delta^{\beta-1}\lara{\xi}}{c_0\lara{\xi}},\qquad\qquad\quad 1+r\le p\le m.
\end{split}
\]
Hence,
\[
|e_{pp}(t,\xi)|\le c'\max\{\eps^{-1},\delta^{\beta-1}\}.
\]
When $p\neq q$ we argue as in \cite{ColKi:02} (estimates (38), (39), (40)). In particular, when both $p$ and $q$ belong to $\{1,...,r\}$ we follow the arguments of Subsection \ref{second_1} for the corresponding term in Case 1. We obtain, for $|\xi|$ large enough,
\[
\begin{split}
|e_{pq}|&\le c\,\delta^{\beta-1}\eps^{\alpha(1-r)},\qquad 1\le p\le m,\, r+1\le q\le m,\\
|e_{pq}|&\le c\,\eps^{\alpha-1},\qquad r+1\le p\le m,\, 1\le q\le r,\\
|e_{pq}|&\le c\,\eps^{-1},\qquad 1\le p\le r,\, 1\le q\le r.
\end{split}
\]
In conclusion,  we get
$$\Vert H^{-1}(t,\xi)\partial_t H(t,\xi)\Vert\le c_2\max\{\eps^{-1},
\delta^{\beta-1}\eps^{\alpha(1-r)}\}$$
for $t\in[0,T]$ and $|\xi|$ large enough.
\subsection{Estimate of the third term}
The matrix $H^{-1}AH$ has entries
\[
c_{pq}(t,\xi)=(\tau_q(t,\xi)-\lambda_q(t,\xi))\frac{\prod_{i=1, i\neq q}^m  (\tau_i(t,\xi)-\lambda_q(t,\xi))}{\prod_{i=1, i\neq p}^m (\lambda_i(t,\xi)-\lambda_p(t,\xi))}.
\]
Arguing as in Case 1, making use of the estimates in Proposition \ref{prop_roots_2} and of the assumption \eqref{hyp_coincide} we obtain, for $|\xi|$ large and $1\le p\le r$, $1\le q\le r$,
\begin{multline}
\label{form_1}
|\tau_q(t,\xi)-\lambda_q(t,\xi)|\frac{\prod_{i=1, i\neq q}^m  |\tau_i(t,\xi)-\lambda_q(t,\xi)|}{\prod_{i=1, i\neq p}^m |\lambda_i(t,\xi)-\lambda_p(t,\xi)|}\\
\le |\tau_q(t,\xi)-\lambda_q(t,\xi)|\frac{\prod_{i=1, i\neq q}^m  |\tau_q(t,\xi)-\lambda_q(t,\xi)|+|\tau_i(t,\xi)-\tau_q(t,\xi)|}{\prod_{i=1, i\neq p}^m |\lambda_i(t,\xi)-\lambda_p(t,\xi)|}\\
=\sum_{k=1}^{m-1}|\tau_q(t,\xi)-\lambda_q(t,\xi)|^k\sum_{S^{(m)}_q(m-k)}\frac{|\tau_{i_1}(t,\xi)-\tau_q(t,\xi)|\cdots
|\tau_{i_{m-k}}(t,\xi)-\tau_q(t,\xi)|}{\prod_{i=1, i\neq p}^m |\lambda_i(t,\xi)-\lambda_p(t,\xi)|}\\
+\frac{|\tau_q(t,\xi)-\lambda_q(t,\xi)|^m}{\prod_{i=1, i\neq p}^m |\lambda_i(t,\xi)-\lambda_p(t,\xi)|}\\
\le c\sum_{k=1}^{m-1}\frac{\eps^{\alpha k}\lara{\xi}^k\lara{\xi}^{m-k}}{\eps^{\alpha(r-1)-\alpha(r-k)}\lara{\xi}^{m-1}}+c\frac{\eps^{\alpha m}\lara{\xi}^m}{\eps^{\alpha(r-1)}\lara{\xi}^{m-1}}\le c'\max\{\eps^\alpha, \eps^{\alpha(m-r+1)}\}\lara{\xi}\\
= c'\eps^\alpha\lara{\xi}.
\end{multline}
If $r+1\le q\le m$ and $1\le p\le r$ then
\begin{multline}
\label{form_2}
|\tau_q(t,\xi)-\lambda_q(t,\xi)|\frac{\prod_{i=1, i\neq q}^m  |\tau_i(t,\xi)-\lambda_q(t,\xi)|}{\prod_{i=1, i\neq p}^m |\lambda_i(t,\xi)-\lambda_p(t,\xi)|}
\le c\delta^\beta\lara{\xi}\frac{1}{\eps^{\alpha(r-1)}} \\
=c\delta^\beta\eps^{\alpha(1-r)}\lara{\xi}.
\end{multline}
If $r+1\le q\le m$ and $1+r\le p\le m$ then
\beq
\label{form_3}
|\tau_q(t,\xi)-\lambda_q(t,\xi)|\frac{\prod_{i=1, i\neq q}^m  |\tau_i(t,\xi)-\lambda_q(t,\xi)|}{\prod_{i=1, i\neq p}^m |\lambda_i(t,\xi)-\lambda_p(t,\xi)|}\le c\delta^\beta\lara{\xi}\frac{1}{c_0}=c'\delta^\beta\lara{\xi}.
\eeq
Finally, if $1\le q\le r$ and $1+r\le p\le m$ then
\beq
\label{form_4}
|\tau_q(t,\xi)-\lambda_q(t,\xi)|\frac{\prod_{i=1, i\neq q}^m  |\tau_i(t,\xi)-\lambda_q(t,\xi)|}{\prod_{i=1, i\neq p}^m |\lambda_i(t,\xi)-\lambda_p(t,\xi)|}\le c\eps^\alpha\lara{\xi}\frac{1}{c_0}=c'\eps^\alpha\lara{\xi}.
\eeq
Combining \eqref{form_1} with \eqref{form_2}, \eqref{form_3} and \eqref{form_4} we obtain
\[
|c_{pq}(t,\xi)|\le c\max\{\eps^\alpha, \delta^\beta\eps^{\alpha(1-r)}, \delta^\beta\}\lara{\xi}=c\max\{\eps^\alpha, \delta^\beta\eps^{\alpha(1-r)}\}\lara{\xi}.
\]
Hence,
$$\Vert H^{-1}AH-(H^{-1}AH)^\ast\Vert\le
c_3\max\{\eps^\alpha, \delta^\beta\eps^{\alpha(1-r)}\}\lara{\xi}$$ for $t\in[0,T]$ and $|\xi|$ large enough.

\subsection{Estimate of the fourth term}
The entries of the matrix $H^{-1}BH$ are given by
\[
d_{pq}(t,\xi)=(-1)^{m-1}\biggl(\prod_{i=1,i\neq p}^m (\lambda_i(t,\xi)-\lambda_p(t,\xi))\biggr)^{-1}g(\lambda_q(t,\xi)),
\]
where $$g(\tau)=\sum_{j=0}^{m-1}(A_{m-j}-A_{(m-j)})(t,\xi)\tau^{j}.$$ Assume that we have lower order terms of order $l$. Then,
\[
\begin{split}
|d_{pq}(t,\xi)|\le c\eps^{\alpha(1-r)}\lara{\xi}^{-m+1+l},\quad & 1\le p\le r,\\
|d_{pq}(t,\xi)|\le c\lara{\xi}^{-m+1+l},\quad & r+1\le p\le m,\\
\end{split}
\]
and $$\Vert H^{-1}BH-(H^{-1}BH)^\ast\Vert\le c_4\eps^{\alpha(1-r)}\lara{\xi}^{l-m+1}$$ for $t\in[0,T]$ and $|\xi|$ large enough.

\subsection{Conclusion of the proof}
We now make use of the four estimates above in \eqref{energy}. We get, for large $|\xi|$,
\begin{multline}
\label{eps_energy_2}
\partial_t |W(t,\xi)|^2\le 2(\rho'(t)\lara{\xi}^{\frac{1}{s}}+c_1\max\{\eps^{-1},\delta^{\beta-1}\}
+c_2\max\{\eps^{-1},\delta^{\beta-1}\eps^{\alpha(1-r)}\} \\ +
c_3\max\{\eps^\alpha, \delta^\beta\eps^{\alpha(1-r)}\}\lara{\xi}+c_4\eps^{\alpha(1-r)}\lara{\xi}^{l-m+1})|W(t,\xi)|^2
\\ +C'\esp^{(\rho(t)-\delta_1)\lara{\xi}^{\frac{1}{s}}}|W(t,\xi)|,
\end{multline}
where $\delta_1>0$ depends on $f$.
Set $\delta=\lara{\xi}^{-1}$ and $\eps=\lara{\xi}^{-\gamma}$. Then we have
\begin{multline}
\label{est_energy_2}
\partial_t |W(t,\xi)|^2\\
\le \left(2\rho'(t)\lara{\xi}^{\frac{1}{s}}
+C\max\{\lara{\xi}^{\gamma},\lara{\xi}^{1-\beta}, \lara{\xi}^{1-\beta-\gamma\alpha(1-r)}, \lara{\xi}^{1-\gamma\alpha}, \lara{\xi}^{-\gamma\alpha(1-r)+l-m+1}\}\right)\cdot\\
\cdot|W(t,\xi)|^2+C'\esp^{(\rho(t)-\delta_1)\lara{\xi}^{\frac{1}{s}}}|W(t,\xi)|\\
=(2\rho'(t)\lara{\xi}^{\frac{1}{s}}+C\max\{\lara{\xi}^{\gamma},\lara{\xi}^{1-\beta-\gamma\alpha(1-r)}, \lara{\xi}^{1-\gamma\alpha}, \lara{\xi}^{-\gamma\alpha(1-r)+l-m+1})\}|W(t,\xi)|^2\\ +
C'\esp^{(\rho(t)-\delta_1)\lara{\xi}^{\frac{1}{s}}}|W(t,\xi)|.
\end{multline}
Let
\[
\gamma=\min\biggl\{\frac{1}{1+\alpha},\frac{\beta}{\alpha r}, \frac{m-l}{\alpha r}\biggr\}.
\]
Hence, $\max\{\gamma, {1-\beta-\gamma\alpha(1-r)}, {-\gamma\alpha(1-r)+l-m+1}\}\le 1-\gamma\alpha$ and
\[
\partial_t |W(t,\xi)|^2\le \left( 2\rho'(t)\lara{\xi}^{\frac{1}{s}}+C\lara{\xi}^{-\gamma\alpha+1}
\right) |W(t,\xi)|^2 + C'\esp^{(\rho(t)-\delta_1)\lara{\xi}^{\frac{1}{s}}}|W(t,\xi)|.
\]
Let $s>0$ be such that
\begin{multline}\label{EQ:formula-s}
\frac{1}{s}>-\min\biggl\{\frac{1}{1+\alpha},\frac{\beta}{\alpha r}, \frac{m-l}{\alpha r}\biggr\}\alpha+1 \\
=\max\biggl\{\frac{1}{1+\alpha}, \frac{r-\beta}{r}, \frac{r-m+l}{r}\biggr\}.
\end{multline}
If $r-m+l>0$, this means that
\beq
\label{s_2a}
s<\min\biggl\{1+\alpha, \frac{r}{r-\beta}, \frac{r}{r-m+l}\biggr\}=1+\min\biggl\{\alpha,\frac{\beta}{r-\beta},\frac{m-l}{r-m+l}\biggr\}.
\eeq
We can assume $|W(t,\xi)|\ge 1$ since when $|W(t,\xi)|<1$ we can use \eqref{VW} to
directly obtain the estimates as in the second line in \eqref{last_estimate2}.
Choosing a suitable decreasing function $\rho$ as in Case 1 we obtain
\beq
\label{energy2}
\partial_t |W(t,\xi)|^2\leq 0
\eeq
for all $t\in[0,T]$ and for $|\xi|$ sufficiently large.  If $r-m+l\le 0$ then the last term under the
maximum sign in \eqref{EQ:formula-s} is negative, and hence disappears. Hence in this case
 \eqref{EQ:formula-s} means that
\beq
\label{s_2b}
s<1+\min\biggl\{\alpha,\frac{\beta}{r-\beta}\biggr\}.
\eeq
Let us finally show that the inequality \eqref{s_2a} is actually also equivalent to
\eqref{s_2b}. Indeed, let us denote $k=m-l$, so that $1\leq k\leq m$.
Consequently, for $\beta\leq 1$ one can readily check that we have
$\frac{\beta}{r-\beta}\leq \frac{k}{r-k}$, proving the claim.

In analogy to Case 1, by arguing as in \eqref{last_estimate}, we see that \eqref{VW} and \eqref{energy2} imply
\beq
\label{last_estimate3}
|V(t,\xi)|\le c\,\esp^{(-\rho(t)+\rho(0))\lara{\xi}^{\frac{1}{s}}}\lara{\xi}^{\gamma\alpha\frac{(r-1)r}{2}}|V(0,\xi)|,
\eeq
for $t\in[0,T]$ and $|\xi|$ large enough.
The estimate \eqref{last_estimate3} proves Theorem \ref{THM:case2}.
Similarly to Case 1, \eqref{last_estimate3} and Proposition \ref{prop_ud}
imply the statement of Theorem \ref{THM:case2u}, also allowing
$s=1+\min\biggl\{\alpha,\frac{\beta}{r-\beta}\biggr\}$.

\begin{rem}\label{REM:r1}
Assume now that we are under assumptions of Case 3, i.e. the Cauchy problem in consideration is strictly hyperbolic. Analysing the estimates of Case 2 under the assumption of strict hyperbolicity, we will set $r=1$ and
repeat the argument first keeping the notation for $\alpha$ and $\beta$
distinguishing them from each other
(although, since we are interested in Case 3,  we will put $\alpha=\beta$ later).
Then, by similar arguments, we readily see that
\begin{enumerate}
\item $|\frac{\partial_t\det H}{\det H}|\le c_1\max\{\eps^{\alpha-1},\delta^{\beta-1}\}$,
\item $\Vert H^{-1}\partial_t H\Vert\le c_2\max\{\eps^{\alpha-1},\delta^{\beta-1}\}$,
\item $\Vert H^{-1}AH-(H^{-1}AH)^\ast\Vert\le c_3\max\{\eps^\alpha,\delta^\beta\}\lara{\xi}$,
\item $\Vert H^{-1}BH-(H^{-1}BH)^\ast\Vert\le c_4\lara{\xi}^{-m+1+l}$,
\end{enumerate}
for $t\in[0,T]$ and $|\xi|$ large enough. Hence, setting $\delta=\lara{\xi}^{-1}$ and $\eps=\lara{\xi}^{-\gamma}$ in the energy estimate
\eqref{eps_energy_2}--\eqref{est_energy_2}
we obtain
\begin{multline*}
\partial_t |W(t,\xi)|^2\le\left(2\rho'(t)\lara{\xi}^{\frac{1}{s}}
+C\max\{\lara{\xi}^{-\gamma\alpha+\gamma},\lara{\xi}^{1-\beta}, \lara{\xi}^{1-\gamma\alpha}, \lara{\xi}^{l-m+1}\}\right)\cdot\\
\cdot|W(t,\xi)|^2+C'\esp^{(\rho(t)-\delta_1)\lara{\xi}^{\frac{1}{s}}}|W(t,\xi)|.
\end{multline*}
Arguing as in Case 2, from $\max\{1-\beta,1-m+l\}\le 1-\gamma\alpha$ we have that $W(t,\xi)$ is of Gevrey order $s$ with
\[
\frac{1}{s}>-\min\biggl\{\frac{\beta}{\alpha}, \frac{m-l}{\alpha}\biggr\}\alpha+1=\max\biggl\{{1-\beta}, {1-m+l}\biggr\}=1-\beta.
\]
This means that $$1\le s<1+\frac{\beta}{1-\beta}.$$
Finally we note that since $m-l\geq 1\geq\beta$, we have in this argument
$\gamma=\min\{\frac{\beta}{\alpha}, \frac{m-l}{\alpha}\}=
\frac{\beta}{\alpha}.$ Recalling that in Case 3, we actually assume $\alpha=\beta$,
we get that in fact $\gamma=1$ (and hence also $\epsilon=\delta$,
simplifying the proof of Case 3 compared to that of Case 2, if needed).
\end{rem}

\end{document}